# A new generalization of the geometric distribution using Azzalini's mechanism: properties and application

(***Preprint**-9$^{th}$ October, 2020*)

*Seng Huat Ong*

Department of Actuarial Science and Applied Statistics, UCSI University, 56000 Kuala Lumpur, Malaysia.

e-mail: ongsh@ucsiuniversity.edu.my

*Subrata Chakraborty*

Department of Statistics, Dibrugarh University, Dibrugarh-786004, Assam, India.

e-mail: subrata_arya@yahoo.co.in

*Aniket Biswas*

Department of Statistics, Dibrugarh University, Dibrugarh-786004, Assam, India.

e-mail: biswasaniket44@gmail.com

**Abstract**

The skewing mechanism of Azzalini for continuous distributions is used **for the first time** to derive a new generalization of the geometric distribution. Various structural properties of the proposed distribution are investigated. Characterizations, including a new result for the geometric distribution, in terms of the proposed model are established. Extensive simulation experiment is done to evaluate performance of the maximum likelihood estimation method. Likelihood ratio test for the necessity of additional skewing parameter is derived and corresponding simulation based power study is also reported. Two real life count datasets are analyzed with the proposed model and compared with some recently introduced two-parameter count models. The findings clearly indicate the superiority of the proposed model over the existing ones in modelling real life count data.

# 1. Introduction

There are various extensions and generalizations of the geometric distribution proposed by many researchers, for instance, Jain and Consul (1971), Philippou et al. (1983), Tripathi et al. (1987), Makcutek (2008), Gómez-Déniz (2010) and Nekoukhou et al. (2012). Many of these generalizations rely on the discretization of continuous distributions or extension using the Marshall-Olkin family of distributions. Recently, Chakraborty and Gupta (2015) proposed a generalization known as the exponentiated geometric distribution and Bhati and Joshi (2018) introduced a weighted geometric distribution.

For continuous distributions many families have been proposed to generate and generalize distributions. In particular, skewed continuous distributions have been thoroughly examined and extended by many researchers stemming from the popular skew normal distribution of Azzalini (1985). The skew normal distribution has probability density function (pdf)

$$f(x) = 2\phi(x)\Phi(\alpha x).$$

As a generalization,

$$f(x) = g(x)G(\alpha x) / \int_{-\infty}^{\infty} g(t)G(\alpha t)dt$$

where $g(x)$ is a pdf with cumulative distribution function (cdf) $G(x)$. This may be viewed as a weighted distribution with the cdf as weighting function. It is of pedagogical interest to consider an analogue of the skewed continuous distributions. To this end, define a probability mass function (pmf) as follows:

$$P_X(x) = p(x)F(\alpha x)/W, \qquad (1.1)$$

where $p(x)$ is a pmf with cdf $F(x)$, $\alpha$ is a non-negative integer and $W$ is the normalizing constant.

*While numerous continuous distributions are introduced in last two decades with Azzalini's set up, no such attempts with discrete distributions has come to our notice so far. The objective of this paper is to define skewed geometric distributions by using skewing mechanism given by (1.1).* In what follows we consider the geometric distribution with pmf

$$p(x) = qp^x, q = 1-p, \ x = 0,1,2,3,... \qquad (1.2)$$

and cdf 
$$F(x) = \sum_{k=0}^{x} p(k) = 1 - p^{x+1}. \qquad (1.3)$$

Formula (1.2) is the probability of *x* successes before the first failure, where the probability of failure is *p* and that of success is *q=1-p*.

Rest of the article is as follows. In the next section, we introduce the proposed distribution and investigate its structural properties. In section 3, we discuss problem of maximum likelihood estimation. Section 4 is about random sample generation and simulation experiment for performance evaluation of MLE. In section 5, we implement likelihood ratio test and conduct a simulation based power study. Two data modelling comparative applications of the proposed distribution is presented in section 6. We conclude the article with a discussion.

## 2. Skewed geometric distribution

On substituting the geometric pmf and cdf into (1.1), we have

$$P_X(x) = qp^x\left(1 - p^{\alpha x+1}\right)/W, \ 0 < p < 1, \ \alpha > 0, \ x \in \{0, 1, 2, \ldots\} \qquad (2.1)$$

where $W = \sum_{k=0}^{\infty} qp^k\left(1 - p^{\alpha k+1}\right) = 1 - \frac{pq}{1 - p^{\alpha+1}} = \frac{1 - p^{\alpha+1} - pq}{1 - p^{\alpha+1}}$.

since $\sum_{k=0}^{\infty} qp^k = 1$ and $\sum_{k=0}^{\infty} pqp^{(\alpha+1)k} = \frac{pq}{1 - p^{(\alpha+1)}}$.

Equation (2.1) is the pmf of a skewed geometric distribution.

**Special cases:**

- If $\alpha = 0$, (2.1) is the geometric pmf.
- When $\alpha \to \infty$

$$P_X(x) = \begin{cases} \dfrac{qp^x}{1 - pq}, & x = 1, 2, \cdots \\ \dfrac{q^2}{1 - pq}, & x = 0. \end{cases}$$

**Remark 1**. In view of the above result it is apparent that for large $\alpha$ the skewed geometric distribution can zero deflate or inflate according as $q \to 0$ or $1$. [See the pmf plots in Figure1 for a quick check].

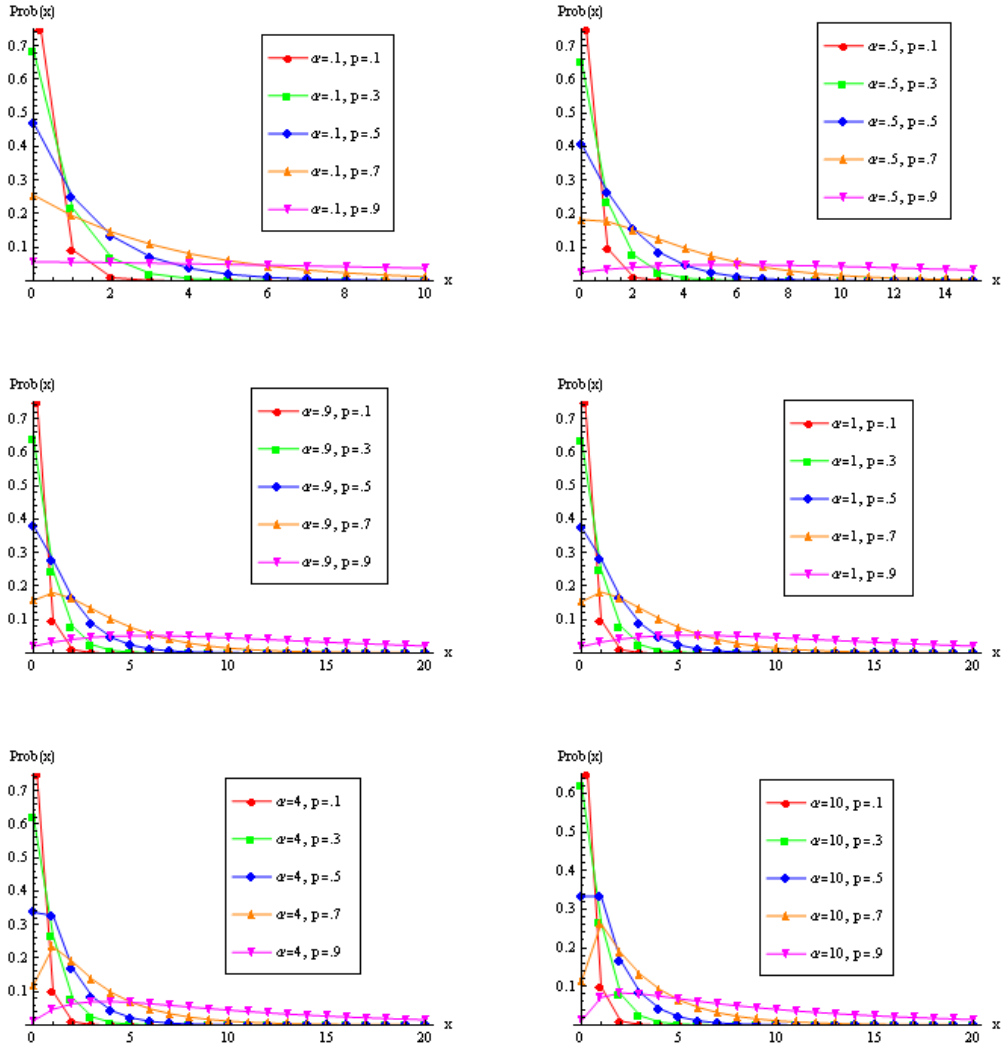

**Figure 1. Plots of probability mass functions**

The skewed geometric pmf (2.1) is still a legitimate pmf if $\alpha$ is a real number. Thus we define

$$P_X(x) = \frac{qp^x\left(1-p^{\alpha x+1}\right)\left(1-p^{(\alpha+1)}\right)}{\left(1-p^{(\alpha+1)}-pq\right)}, \ 0<p<1, \ \alpha>0, \ x=0,1,2,3,... \quad (2.2)$$

as the skewed geometric pmf and we denote it by $SG(p,\alpha)$.

**Remark 2.** By using a weighted distribution construction, Bhati and Joshi (2018) considered a weighting function $\omega(x) = 1 - p^{\alpha(x+1)}$ for the geometric distribution. Our approach based on (1.1) uses the cdf as the weighting function results in a different weight $\omega(x) = 1 - p^{\alpha x+1}$. Thus (2.2) gives a generalized geometric distribution which differs from the one proposed by these authors

As a generalization of (2.2), replace the geometric pmf with an arbitrary pmf $p(x)$. This results in

$$P_X(x) = p(x)(1 - p^{\alpha x+1})/W \qquad (2.3)$$

where $W = \sum_{k=0}^{\infty} p(k)(1 - p^{\alpha k+1})$.

Let $H(z) = \sum_{x=0}^{\infty} p(x) z^x$ be the probability generating function (pgf) for $p(x)$. Then the pgf of (2.3) is given by

$$G_X(z) = \sum_{x=0}^{\infty} P_X(x) z^x = \left[H(z) - pH(p^{\alpha} z)\right]/W, \qquad (2.4)$$

and $W = \sum_{k=0}^{\infty} p(k)(1 - p^{\alpha k+1}) = \sum_{k=0}^{\infty} p(k) - \sum_{k=0}^{\infty} p(k) p^{\alpha k+1} = 1 - pH(p^{\alpha})$.

**Probability recurrence relation and tail length:**

Consider the ratio of probabilities

$$\frac{P_X(x+1)}{P_X(x)} = \frac{p(1 - p^{\alpha x + \alpha + 1})}{(1 - p^{\alpha x + 1})}.$$

This gives a simple two-term recurrence formula

$$P_X(x+1) = p\{(1 - p^{\alpha x + \alpha + 1})/(1 - p^{\alpha x + 1})\} P_X(x).$$

The behaviour of the tail of the distribution may be deduced from the ratio of probabilities (Ong and Muthaloo, 1995). As $x \to \infty$ the above tends to $p \neq 0$. Hence the distribution is short/long tailed depending on the value of $p$. As such when $p \to 0$ we get a Poisson-type tail whereas for $p \to 1$, a longer tail appears.

### 2.1 Characterisation results

A few stochastic characterisations of the skewed geometric distribution are investigated here.

**Theorem1.** Let $X_1, X_2$ be independent and identically distributed (i.i.d.) discrete random variables, For any positive integer $\alpha$ the conditional distribution $P(X_1 = x | X_2 \leq \alpha x)$ is skewed geometric $SG(p, \alpha)$ with pmf (2.2) if, and only if, $X_1, X_2$ are geometric random variables defined by (1.2).

**Proof:**

*If part*: Let $X_1, X_2$ be i.i.d. Geometric $(p)$. Then

$$P(X_1 = x | X_2 \leq \alpha X_1) = \frac{P(X_1 = x \cap X_2 \leq \alpha X_1)}{P(X_2 \leq \alpha X_1)}$$

$$= \frac{P(X_1 = x)P(X_2 \leq \alpha x)}{P(X_2 \leq \alpha X_1)} = \frac{q p^x (1 - p^{\alpha x + 1})}{P(X_2 \leq \alpha X_1)} \quad (2.5)$$

But $P(X_2 \leq \alpha X_1) = \sum_{j \geq 0} P(X_2 \leq \alpha j) P(X_1 = j) = W$

Therefore (2.5) becomes

$$P(X_1 = x | X_2 \leq \alpha X_1) = \frac{q p^x (1 - p^{\alpha x + 1})(1 - p^{\alpha + 1})}{(1 - p^{\alpha + 1} - pq)}.$$

*Only if part*:

Again $P(X_1 = x | X_2 \leq \alpha X_1) = \frac{q p^x (1 - p^{\alpha x + 1})(1 - p^{\alpha + 1})}{(1 - p^{\alpha + 1} - pq)}.$

We rearrange the above as

$$P(X_1 = x) F_{X_2}(\alpha x) = P(X_2 \leq \alpha X_1) \frac{q p^x (1 - p^{\alpha x + 1})(1 - p^{\alpha + 1})}{(1 - p^{\alpha + 1} - pq)}$$

In particular, for $\alpha = 1$

$$P(X_1 = x) F_{X_2}(x) = P(X_2 \leq X_1) \frac{q p^x (1 - p^{x+1})(1 - p^2)}{(1 - p^2 - pq)} = C p^x (1 - p^{x+1}), say$$

when $x = 0$, 1

$P(X_1 = 0) = C(1 - p)$

$P(X_1 = 1)\{P(X_1 = 0) + P(X_1 = 1)\} = C p(1 - p^2)$

$\Rightarrow \theta_1 P(X_1 = 0) P(X_1 = 0) + \theta_1^2 P(X_1 = 0)^2 = C p(1 - p^2),$

where $\theta_1 = \frac{P(X_1 = 1)}{P(X_1 = 0)}.$

Solving these two simultaneously for $\theta_1$ we get $\theta_1 = p$ or $-1 - p$.

Hence $P(X_1 = 1) = p P(X_1 = 0)$.

We claim that $P(X_1 = y) = p^y P(X_1 = 0)$. That this is true can be seen from last equation.

We shall now proceed through mathematical induction assuming that for a positive integer $n$

$$P(X_1 = n) = p^n P(X_1 = 0) \quad (2.6)$$

Now

$$P(X_1 = n+1) F_{X_2}(n+1) = C p^{n+1}(1 - p^{n+2})$$

$$P(X_1 = n+1)\{P(X_1 = 0) + P(X_1 = 1) + \cdots + P(X_1 = n+1)\} = C\, p^{n+1}(1-p^{n+2})$$

which reduces to

$$\theta_{n+1} P(X_1 = 0)\{P(X_1 = 0) + pP(X_1 = 0) + \cdots + p^n P(X_1 = 0)\} + \theta_{n+1} P(X_1 = 0) = C\, p^{n+1}(1-p^{n+2})$$

On taking, $\theta_k = P(X_1 = k)/P(X_1 = 0)$. (2.7)

Solving this we get $\theta_{n+1} = p^{n+1}$ or $-\sum_{j=0}^{n+1} p^j$.

Thus from (2.7)

$$\theta_{n+1} = \frac{P(X_1 = n)}{P(X_1 = 0)} \Rightarrow P(X_1 = n) = p^{n+1} P(X_1 = 0).$$

Hence $P(X_1 = y) = (1-p) p^y \sim$ Geometric $(p)$.

**Theorem 2.** Let $X$ be any discrete random variable with support $\mathbf{N} \cup \{0\}$ then $X$ follows skewed geometric distribution $SG(p, \alpha)$ in (2.2) if, and only if,

$$E\left[\frac{p^x(1+p)(1-p^{\alpha+1}) - p^{x(\alpha+1)+1}(1-p)(1+p^{\alpha+1})}{1-p+p^2-p^{\alpha+1}} \,\bigg|\, X \geq k\right] = \frac{p^k(1-p^{\alpha+1} - p^{\alpha k+1} + p^{\alpha k+2})}{(1-p+p^2-p^{\alpha+1})} \quad \text{(A)}$$

**Proof:**

If $X$ follows $SG(p, \alpha)$ then left-hand side of (A) can be written as

$$\sum_{x \geq k} \frac{p^x(1+p)(1-p^{\alpha+1}) - p^{x(\alpha+1)+1}(1-p)(1+p^{\alpha+1})}{1-p+p^2-p^{\alpha+1}} \cdot P[X = x \mid X \geq k]$$

$$= \sum_{x \geq k} \frac{p^x(1+p)(1-p^{\alpha+1}) - p^{x(\alpha+1)+1}(1-p)(1+p^{\alpha+1})}{1-p+p^2-p^{\alpha+1}} \cdot \frac{p^x(1-p)(1-p^{\alpha x+1})(1-p^{\alpha+1})}{(1-p^{\alpha+1} - pq)} \bigg/ P[X \geq k]$$

$$= \sum_{x \geq k} \frac{p^x(1+p)(1-p^{\alpha+1}) - p^{x(\alpha+1)+1}(1-p)(1+p^{\alpha+1})}{1-p+p^2-p^{\alpha+1}} \cdot \frac{p^x(1-p)(1-p^{\alpha x+1})(1-p^{\alpha+1})}{p^k(1-p^{\alpha+1} - p^{\alpha k+1} + p^{\alpha k+2})}$$

After some algebra we get

$$= \frac{p^{2k}(1-p^{\alpha+1} - p^{\alpha k+1} + p^{\alpha k+2})^2}{(1-p+p^2-p^{\alpha+1})^2} \cdot \frac{1-p+p^2-p^{\alpha+1}}{p^k(1-p^{\alpha+1} - p^{\alpha k+1} + p^{\alpha k+2})}$$

$$= \frac{p^k(1-p^{\alpha+1} - p^{\alpha k+1} + p^{\alpha k+2})}{1-p+p^2-p^{\alpha+1}}.$$

Conversely, assuming (A) to be correct

$$\frac{1}{P(X \geq k)} \sum_{x \geq k}^{\infty} \frac{p^x(1+p)(1-p^{\alpha+1}) - p^{x(\alpha+1)+1}(1-p)(1+p^{\alpha+1})}{1-p+p^2-p^{\alpha+1}} \cdot P[X=x]$$

$$= \frac{p^k(1-p^{\alpha+1} - p^{\alpha k+1} + p^{\alpha k+2})}{1-p+p^2-p^{\alpha+1}}$$

Therefore,

$$\frac{1}{P(X \geq k)} \sum_{x \geq k+1}^{\infty} \frac{p^x(1+p)(1-p^{\alpha+1}) - p^{x(\alpha+1)+1}(1-p)(1+p^{\alpha+1})}{1-p+p^2-p^{\alpha+1}} P[X=x]$$

$$= \frac{p^{k+1}(1-p^{\alpha+1} - p^{\alpha k+1} + p^{\alpha(k+1)+2})}{1-p+p^2-p^{\alpha+1}}$$

Taking difference of the preceding two equations we get

$$\frac{p^k(1-p^{\alpha+1} - p^{\alpha k+1} + p^{\alpha k+2})}{1-p+p^2-p^{\alpha+1}} \cdot \frac{P[X=k]}{P(X \geq k)} = \frac{p^k(1-p)(1-p^{\alpha+1})(1-p^{\alpha k+1})}{1-p+p^2-p^{\alpha+1}}$$

$$\frac{P[X=k]}{P(X \geq k)} = \frac{(1-p)(1-p^{\alpha+1})(1-p^{\alpha k+1})}{(1-p^{\alpha+1} - p^{\alpha k+1} + p^{\alpha k+2})}.$$

This is $SG(p, \alpha)$. Hence the result is proved.

**Theorem 3.** Let $X$ be any discrete random variable with support $\mathbf{N} \cup \{0\}$ then $X$ follows $SG(p, \alpha)$ in in (2.2) if, and only if,

$$h_X(x+1) - h_X(x) = \frac{(1-p)(1-p^{\alpha+1})(p^{x\alpha} - 2p^{(x+1)\alpha} + p^{(x+2)\alpha})}{(1-p^{\alpha+1} - p^{(x+1)\alpha+1} + p^{(x+1)\alpha+2})(1-p^{\alpha+1} - p^{(x+2)\alpha+1} + p^{(x+2)\alpha+2})}$$

with initial condition $h_X(0) = \dfrac{(1-p)^2(1-p^{\alpha+1})}{p(1-2p^{\alpha+1} - p^{\alpha+2})}$.

**Proof:** When $X$ follows $SG(p, \alpha)$ in (2.2) the result can be easily seen to be true.

Conversely assuming the result to be true we can write

$$\sum_{x=0}^{k} \{h_X(x+1) - h_X(x)\} = \sum_{x=0}^{k} \frac{(1-p)(1-p^{\alpha+1})(p^{x\alpha} - 2p^{(x+1)\alpha} + p^{(x+2)\alpha})}{(1-p^{\alpha+1} - p^{(x+1)\alpha+1} + p^{(x+1)\alpha+2})(1-p^{\alpha+1} - p^{(x+2)\alpha+1} + p^{(x+2)\alpha+2})}$$

that is,

$$\sum_{x=0}^{k} \{h_X(x+1) - h_X(x)\} =$$

$$\sum_{x=0}^{k} \left\{ \frac{(1-p)(1-p^{\alpha+1})(1-p^{\alpha(x+1)+1})}{p(1-p^{\alpha+1} - p^{\alpha(x+2)+1} + p^{\alpha(x+2)+2})} - \frac{(1-p)(1-p^{\alpha+1})(1-p^{\alpha x+1})}{p(1-p^{\alpha+1} - p^{\alpha(x+1)+1} + p^{\alpha(x+1)+2})} \right\},$$

$$h_X(k+1) - h_X(0) = \frac{(1-p)(1-p^{\alpha+1})(1-p^{\alpha(k+1)+1})}{p(1-p^{\alpha+1}-p^{\alpha(k+2)+2})} - \frac{(1-p)(1-p^{\alpha+1})(1-p)}{p(1-p^{\alpha+1}-p^{\alpha+1}+p^{\alpha+2})}$$

This implies

$$h_X(k+1) = \frac{(1-p)(1-p^{\alpha+1})(1-p^{\alpha(k+1)+1})}{p(1-p^{\alpha+1}-p^{\alpha(k+2)+2})} - \frac{(1-p)(1-p^{\alpha+1})(1-p)}{p(1-p^{\alpha+1}-p^{\alpha+1}+p^{\alpha+2})} + \frac{(1-p)^2(1-p^{\alpha+1})}{p(1-2p^{\alpha+1}-p^{\alpha+2})}$$

On simplification we get $h_X(k+1) = \dfrac{(1-p)(1-p^{\alpha+1})(1-p^{\alpha(k+1)+1})}{p(1-p^{\alpha+1}-p^{\alpha(k+2)+2})}$.

This is the hazard function for $SG(p, \alpha)$.

### 2.2 Cumulative distribution function

By using the geometric cdf $F(x) = 1 - p^{x+1}$, it is easy to work out the formula for cdf of the $SG(p, \alpha)$ as

$$F_X(x) = \sum_{j=0}^{x} qp^j \left(1 - p^{\alpha j+1}\right) / W, \quad W = 1 - \frac{pq}{1-p^{\alpha+1}}$$

But we know that

$$\sum_{j=0}^{x} qp^j = 1 - p^{x+1} \text{ and } \sum_{j=0}^{x} pqp^j = pq \sum_{j=0}^{x} p^{(\alpha+1)j} / (1 - p^{\alpha+1}) = pq(1 - p^{(\alpha+1)(x+1)}) / (1 - p^{\alpha+1}).$$

And this gives,

$$F_X(x) = \sum_{j=0}^{x} qp^j (1 - p^{\alpha j+1}) / W = \left\{(1 - p^{x+1}) - pq \frac{1 - p^{(\alpha+1)(x+1)}}{1 - p^{\alpha+1}}\right\} / W$$

### 2.3 Probability generating function

The pgf of the $SG(p, \alpha)$ is given by

$$G_X(z) = \left[H(z) - pH(p^\alpha z)\right] / W, \quad W = 1 - \frac{pq}{1-p^{\alpha+1}} \tag{2.8}$$

where $H(z) = \left(\dfrac{q}{1-pz}\right)$ is the geometric pgf.

### 2.4 Moments

By differentiating (2.8) $r$ times with respect to $z$, and setting $z = 1$, the $r^{\text{th}}$ factorial moment can be obtained.

The mean of $SG(p, \alpha)$ is the first factorial moment.

$H'(z)\big|_{z=1} = \dfrac{p}{1-p}$, is the mean of geometric distribution in (1.2).

$$H'(p^\alpha z)\big|_{z=1} = \sum_{k=0}^{\infty} kp(k)p^{\alpha k}z^{k-1}\bigg|_{z=1} = \sum_{k=0}^{\infty} kp(k)p^{\alpha k} = \frac{q}{(1-p^{\alpha+1})}\sum_{k=0}^{\infty} k(1-p^{\alpha+1})p^{(\alpha+1)k} = \frac{qp^{\alpha+1}}{(1-p^{\alpha+1})^2}$$

**Mean of** $SG(p,\alpha)$:

$$\mu_X = \left[\frac{p}{1-p} - pqp^{\alpha+1}(1-p^{\alpha+1})^{-2}\right]\bigg/\left(1-\frac{pq}{1-p^{\alpha+1}}\right) = \frac{p[1-3p^{\alpha+1}+2p^{\alpha+2}-p^{\alpha+3}+p^{2\alpha+2}]}{(1-p)(1-p^{\alpha+1})(1-p+p^2-p^{\alpha+1})}$$

and $Var(X) = \left[-\frac{2p^2}{(1-p)^2} - \frac{2qp^{3+2\alpha}}{(1-p^{1+\alpha})^3}\right]\bigg/\left(1-\frac{pq}{1-p^{\alpha+1}}\right) - \mu_X^2$

Similarly, one can derive expression for variance and other higher order moments. The distribution as expected is over dispersed and index of dispersion can be very high for small vales of $\alpha$ and higher values of $p$.

### 2.5 Log-concavity and unimodality

A distribution is said to be log-concave if its pmf $\{f_k\}$, $f_k > 0$, $\forall k$ satisfies $f_k^2 \geq f_{k+1}f_{k-1}$.
The skewed geometric distribution is log-concave because

$$\{qp^x(1-p^{\alpha x+1})\}^2 \geq qp^{x+1}(1-p^{\alpha(x+1)+1})qp^{x-1}(1-p^{\alpha(x-1)+1})$$

This follows from $(p^\alpha - 1)^2 \geq 0$, $0 < p < 1, \alpha > 0$.

Theorem 3 of Keilson and Geber (1971, page 386) states that a necessary and sufficient condition for pmf $\{f_k\}$ be strongly unimodal is that $f_k$ be log-concave for all $k$. Thus the skewed geometric distribution is strongly unimodal.

**Remark3**. It can be verified that the distribution has a nonzero mode if $\alpha > \log(2p-1)/\log(p) - 1$. For example when $p = .7$, the distribution will have nonzero mode for $\alpha > 0.56898$. This can be observed in the last four frames in the Figure 1.

### 2.6 Reliability properties

The survival and failure rate functions of the $SG(p,\beta)$ are respectively given by

$$\bar{F}_X(x) = 1 - \sum_{j=0}^{x} qp^j(1-p^{\alpha j+1})/W = 1 - \left\{(1-p^{x+1}) - pq\frac{1-p^{(\alpha+1)(x+1)}}{1-p^{\alpha+1}}\right\}\bigg/W$$

and $h_X(x) = \frac{P_X(x)}{\bar{F}_X(x)} = \frac{q}{p}\frac{(1-p^{\alpha x+1})(1-p^{(\alpha+1)})}{(1-p^{\alpha+1}-qp^{\alpha(x+1)+1})}$.

The failure rate function plotted in Figure 2 for some chosen values of parameters reveals its monotonic increasing nature.

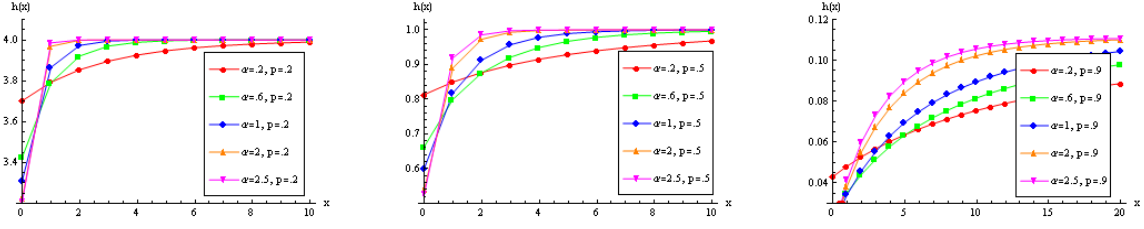

**Figure 2. Plots of Failure rate function**

Moreover it follows from log-concavity that $SG(p,\alpha)$ has an increasing failure rate (*IFR*) (Gupta et al., 2008). Furthermore we have

$$IFR \Rightarrow IFRA \Rightarrow NBU \Rightarrow NBUE \Rightarrow HNBUE$$

where *IFRA* (increasing failure rate average), *NBU* (new better than used), *NBUE* (new better than used in expectation) and *HNBUE* (harmonic new better than used in expectation). Therefore, the skewed geometric distribution is *IFR*, *IFRA*, *NBU*, *NBUE* and *HNBUE*.

## 3. Parameter estimation

In this section, maximum likelihood estimation is considered for parameter estimation. For a random sample $\boldsymbol{x} = (x_1, x_2, \ldots, x_n)$ of size $n$ from $SG(\alpha, p)$ the log-likelihood function is given by

$$P_X(x) = \frac{qp^x\left(1-p^{\alpha x+1}\right)\left(1-p^{\alpha+1}\right)}{\left(1-p^{\alpha+1}-pq\right)}, \quad 0<p<1, \ \alpha>0, \ x=0,1,2,3,\ldots$$

$$l(p,\alpha|\boldsymbol{x}) = \log(p,\alpha|\boldsymbol{x})$$
$$= n\log(1-p) + \log(p)\sum_{i=1}^{n} x_i + \sum_{i=1}^{n} \log(1-p^{\alpha x_i+1}) + n\log(1-p^{\alpha+1}) - n\log(1-p^{\alpha+1}-p(1-p))$$

The corresponding score functions are

$$\frac{\partial l(p,\alpha|\boldsymbol{x})}{\partial p} = -\frac{n}{1-p} + \frac{1}{p}\sum_{i=1}^{n} x_i - \sum_{i=1}^{n} \frac{(\alpha x_i+1)p^{\alpha x_i}}{(1-p^{\alpha x_i+1})} - \frac{n(\alpha+1)p^\alpha}{(1-p^{\alpha+1})} - \frac{n[(\alpha+1)p^\alpha - 2p+1]}{1-p^{\alpha+1}-p(1-p)} \quad \text{and}$$

$$\frac{\partial l(p,\alpha|\boldsymbol{x})}{\partial \alpha} = -\sum_{i=1}^{n} \frac{x_i p^{\alpha x_i+1}\log p}{1-p^{\alpha x_i+1}} - \frac{np^{\alpha+1}\log p}{(1-p^{\alpha+1})} + \frac{np^{\alpha+1}\log p}{1-p^{\alpha+1}-p(1-p)}$$

Solving these ML equations simultaneously for $(p,\alpha)$ is not feasible. This directed us to try numerical maximization of $L(p,\alpha|\boldsymbol{x})$ with respect to p and α. But due to the lack of suitable initial values for the parameters to be estimated, numerical routines often produced inconsistent outputs. Exhaustive search over the parameter space of $(p,\alpha)$ in such a scenario a reasonable option, but will be a futile exercise here due to the unbounded nature of the parameter space which is $(0, 1) \times [0, \infty)$. We address this issue in the next subsection.

### 3.1 A useful re-parameterisation

In the pmf stated in (2.2) the parameter $\alpha$ in not bounded, this causes problem in setting an initial value for successfully running a maximization routine to obtain MLEs. In order to alleviate this issue we consider a simple re-parameterisation by substituting $p^\alpha$ by a new parameter $\beta$. The resulting pmf is then given as

$$P_X(x) = \frac{p^x(1-p)(1-p\beta)(1-p\beta^x)}{1-p\beta-p(1-p)}, \quad x = 0,1,\dots\,; 0 < p < 1, 0 < \beta < 1. \tag{2.9}$$

This re-parameterization will not alter any result presented in the preceding sections and can be easily expressed in terms of $(p, \beta)$ instead of $(\alpha, p)$. We refer to the pmf in (2.9) as $RSG(p, \beta)$.

Note that $0 < \beta < 1$, means we can employ exhaustive search algorithm without bothering about initial starting value as in case of numerical optimization. For a sample $(x_1, x_2, \dots, x_n)$ of size $n$ from $RSG(p, \beta)$. The log likelihood function is given by

$$l(p, \beta | x) = log(p, \beta | x)$$

$$= n \log(1-p) + \log(p) \sum_{i=1}^{n} x_i + \sum_{i=1}^{n} \log(1 - p\beta^{x_i})$$
$$+ n \log(1 - p\beta) - n \log(1 - p\beta - p(1-p))$$

Given $x$, the parameters are estimated (correct up to 4 decimal places) using the following search algorithm:

(i) Fix $p \in \mathbf{p} = \{0.0001, 0.0002, \dots, 0.9999\}$
(ii) For each $p$ in $\mathbf{p}$, find $\beta \in \boldsymbol{\beta} = \{0.0001, 0.0002, \dots, 0.9999, 1.000\} \ni L(p, \beta | x)$ is maximum.
(iii) $\forall$ pairs of $(p, \beta) \in (\mathbf{p}, \boldsymbol{\beta})$ obtained via step (ii), find the pair $(\hat{p}, \hat{\beta})$ for which $L(p, \beta | x)$ is maximum.

For using asymptotic properties of the MLEs $(\hat{p}, \hat{\beta})$ we need the following:

$$\frac{\partial^2 L(p, \beta | x)}{\partial p^2} = -\frac{n}{(1-p)^2} - \frac{n\beta^2}{(1-\beta p)^2} + \frac{n(\beta^2 - 2p(1-p) + 2\beta(1-p) - 1)}{(1-(1+\beta)p + p^2)^2}$$
$$- \sum_{i=1}^{n} \frac{\beta^{2x_i}}{(1-p\beta^{x_i})^2} - \frac{1}{p^2} \sum_{i=1}^{n} x_i = J_{11}$$

$$\frac{\partial^2 L(p,\beta|x)}{\partial p \partial \beta} = -\frac{np(p(4 + \beta(2 - 4p) + (-2 + p)p + \beta^2 p^2) - 2)}{(1 - \beta p)^2((1 + \beta - p)p - 1)^2} - \sum_{i=1}^{n}\left(\frac{px_i \beta^{2x_i - 1}}{(1 - p\beta^{x_i})^2}\right.$$

$$\left. + \frac{x_i \beta^{x_i - 1}}{1 - p\beta^{x_i}}\right) = J_{12} = J_{21}$$

$$\frac{\partial^2 L(p,\beta|x)}{\partial \beta^2} = n\beta^2 \left(\frac{1}{(1 - (1 + \beta)p + p^2)^2} - \frac{1}{(1 - \beta p)^2}\right)$$

$$- \sum_{i=1}^{n} \frac{px_i(x_i - 1)\beta^{x_i - 2}}{1 - p\beta^{x_i}} + \frac{p^2 x_i^2 \beta^{2x_i - 2}}{(1 - p\beta^{x_i})^2} = J_{22}$$

Writing $J(p,\beta|x) = \begin{pmatrix} J_{11} & J_{12} \\ J_{21} & J_{22} \end{pmatrix}$ an approximation of the Fisher's information matrix is obtained as $I = -J(\hat{p},\hat{\beta}|x)$. An approximate dispersion matrix of the MLEs $(\hat{p},\hat{\beta})$ is then computed by $\Sigma = I^{-1}$ from which we obtain the variances and covariance of the parameter estimates.

## 4. Simulation based performance study of MLE
## 4.1 Random sample generation and

As an explicit from of the cdf is available, one can use probability integral transform to generate a random observations $x$ from $SG(p,\alpha)$ by solving $F_X(x) = u$ for $x$, given the values of $p, \alpha$ and $U = u \sim Uniform(0,1)$. In fact, $x = F_X^{-1}(u)$, where $F_X^{-1}(.)$ is of course the corresponding quantile function. But, for the proposed distribution a direct solution for $x$ as above is not feasible and the various numerical routines lack consistency and are found to be imprecise. In view of this situation, we prescribe the following exhaustive search based algorithm that provides consistently precise output.

*Algorithm*

(For generating a sample of size $n$ from $SG(p, \alpha)$)

i. Input $p, \alpha$

ii. Compute : $W = W(p,\alpha)$, from equation (2.3), $C = \dfrac{p(1-p)}{1 - p^{\alpha+1}}$

iii. Generate one $u \sim Uniform(0,1)$ and Compute $B = 1 - C - uW$

iv. For $y = \{0.0001, 0.0002, ..., 0.9999\}$, obtain $y_0 = \min_{y}\{y - C.y^{\alpha+1} - B\}$

v. Find $x = \left\lfloor \dfrac{\log(y_0)}{\log(p)} \right\rfloor$

vi. Repeat steps (i) to (v) $n$ times.

Obviously, generating a random sample from $RSG(p, \beta)$ is equivalent to generating the same from $SG(p, \log \beta / \log p)$.

## 4.2. Simulation experiment

We conduct simulation to study the performance of the MLE by generating random samples of size 50, 100, 200, 300 from $RSG(p, \beta)$. and estimating by ML method. Average bias, MSE, CI are computed by replicating the experiment 1000 times. The findings are presented in Table 1.

**Table 1**

| Sample size | Bias($\hat{p}$) | MSE($\hat{p}$) | CI($\hat{p}$) | Bias($\hat{\beta}$) | MSE($\hat{\beta}$) | CI($\hat{\beta}$) |
|---|---|---|---|---|---|---|
| \multicolumn{7}{c}{$p = 0.2, \beta = 0.2$} |
| 50 | 0.0076 | 0.0020 | (0.0671, 0.3481) | 0.2004 | 0.2333 | (-3.5686, 4.3695) |
| 100 | 0.0071 | 0.0013 | (0.1063, 0.3080) | 0.1743 | 0.2162 | (-2.1264, 2.8749) |
| 200 | 0.0087 | 0.0006 | (0.1172, 0.3002) | 0.1888 | 0.2612 | (-2.1558, 2.9335) |
| 300 | 0.0077 | 0.0005 | (0.1385, 0.2769) | 0.1864 | 0.2160 | (-1.6828, 2.4555) |
| 400 | 0.0072 | 0.0005 | (0.1400, 0.2744) | 0.2035 | 0.2142 | (-1.6186, 2.4255) |
| 500 | 0.0068 | 0.0004 | (0.1481, 0.2654) | 0.1643 | 0.2004 | (-1.3704, 2.0989) |
| \multicolumn{7}{c}{$p = 0.2, \beta = 0.5$} |
| 50 | -0.0016 | 0.0022 | (0.0576, 0.3392) | -0.1014 | 0.2128 | (-3.5133, 4.3105) |
| 100 | 0.0036 | 0.0011 | (0.0997, 0.3076) | -0.0464 | 0.2038 | (-2.1988, 3.1059) |
| 200 | 0.0007 | 0.0007 | (0.1163, 0.2851) | -0.0650 | 0.1983 | (-2.1116, 2.9816) |
| 300 | 0.0054 | 0.0005 | (0.1318, 0.2791) | -0.0040 | 0.1794 | (-1.5181, 2.5100) |
| 400 | 0.0030 | 0.0004 | (0.1307, 0.2753) | -0.0151 | 0.1761 | (-1.5029, 2.4726) |
| 500 | 0.0027 | 0.0003 | (0.1371, 0.2682) | -0.0164 | 0.1867 | (-1.3050, 2.2721) |
| \multicolumn{7}{c}{$p = 0.2, \beta = 0.8$} |
| 50 | -0.0053 | 0.0023 | (0.0051, 0.0043) | -0.3514 | 0.3332 | (-3.3991, 4.2963) |
| 100 | -0.0085 | 0.0012 | (0.0864, .2967) | -0.2870 | 0.2899 | (-2.2936, 3.3194) |
| 200 | -0.0039 | 0.0007 | (0.0983, 0.2940) | -0.2802 | 0.2706 | (-0.2117, 3.1514) |
| 300 | -0.0030 | 0.0005 | (0.1176, 0.2763) | -0.2423 | 0.2537 | (-1.6164, 2.7317) |
| 400 | -0.0020 | 0.0004 | (0.1141, 0.2820) | -0.1941 | 0.2186 | (-1.5557, 2.7675) |
| 500 | -0.0038 | 0.0003 | (0.1254, 0.2670) | -0.2242 | 0.2345 | (-1.3207, 2.4723) |
| \multicolumn{7}{c}{$p = 0.5, \beta = 0.2$} |
| 50 | 0.0068 | 0.0016 | (0.4122, 0.6015) | 0.1061 | 0.1138 | (-0.6873, 1.2996) |
| 100 | 0.0015 | 0.0008 | (0.4419, 0.5611) | 0.0666 | 0.0779 | (-0.4902, 1.0233) |
| 200 | 0.0034 | 0.0004 | (0.4635,0.5433) | 0.0371 | 0.0500 | (-0.3281,0.8031) |
| 300 | 0.0004 | 0.0003 | (0.4693,0.5314) | 0.0103 | 0.0328 | (-0.2605,0.6811) |
| 400 | 0.0004 | 0.0002 | (0.4734,0.5273) | 0.0126 | 0.0301 | (-0.1894,0.6145) |
| 500 | 0.0006 | 0.0001 | (0.4766,0.5246) | 0.0047 | 0.0266 | (-0.1579,0.5672) |
| $p=0.5, \beta=0.5$ | | | | | | |
| 50 | 0.0094 | 0.0021 | (0.3950,0.6293) | -0.0380 | 0.1238 | (-0.4585,1.3824) |
| 100 | 0.0069 | 0.0011 | (0.4263,0.5875) | -0.0377 | 0.0998 | (-0.2927,1.2174) |
| 200 | 0.0036 | 0.0005 | (0.4582,0.5490) | -0.0267 | 0.0582 | (-0.0377,0.9844) |
| 300 | 0.0029 | 0.0004 | (0.4650,0.5409) | -0.0219 | 0.0490 | (0.0567,0.8995) |
| 400 | 0.0018 | 0.0002 | (0.4722,0.5314) | -0.0229 | 0.0365 | (0.1149,0.8393) |
| 500 | 0.0031 | 0.0002 | (0.4753,0.5273) | -0.0132 | 0.0263 | (0.1665,0.8071) |
| $p=0.5, \beta=0.8$ | | | | | | |
| 50 | 0.0074 | 0.0025 | (0.3464,0.6684) | -0.1216 | 0.1242 | (-0.1273,1.4841) |
| 100 | 0.0067 | 0.0017 | (0.3907,0.6227) | -0.0998 | 0.0933 | (0.0811,1.3193) |
| 200 | 0.0079 | 0.0010 | (0.4063,0.6094) | -0.0448 | 0.0473 | (0.2962,1.2141) |
| 300 | 0.0074 | 0.0008 | (0.4297,0.5851) | -0.0277 | 0.0234 | (0.4068,1.1378) |
| 400 | 0.0062 | 0.0007 | (0.4396,0.5728) | -0.0125 | 0.0246 | (0.4685,1.1064) |
| 500 | 0.0061 | 0.0006 | (0.4489,0.5633) | -0.0156 | 0.0218 | (0.5029,1.0659) |
| $p=0.8, \beta=0.2$ | | | | | | |
| 50 | -0.0063 | 0.0007 | (0.7422,0.8451) | 0.0971 | 0.0979 | (-0.4496,1.0439) |
| 100 | -0.0027 | 0.0003 | (0.7605,0.8342) | 0.0641 | 0.0613 | (-0.3473,0.8755) |
| 200 | -0.0022 | 0.0002 | (0.7716,0.8239) | 0.0164 | 0.0365 | (-0.2378,0.6705) |
| 300 | -0.0006 | 0.0001 | (0.7783,0.8206) | 0.0036 | 0.0286 | (-0.1632,0.5705) |
| 400 | -0.0007 | 0.0001 | (0.7809,0.8178) | 0.0006 | 0.0212 | (-0.1191,0.5204) |
| 500 | -0.0010 | 0.0001 | (0.7826,0.8155) | 0.0119 | 0.0189 | (-0.0722,0.4960) |
| $p=0.8, \beta=0.5$ | | | | | | |
| 50 | -0.0001 | 0.0006 | (0.7505,0.8493) | -0.0302 | 0.0880 | (-0.1685,1.1081) |
| 100 | -0.0023 | 0.0003 | (0.7621,0.8332) | -0.0313 | 0.0591 | (-0.0051,0.9425) |
| 200 | 0.0002 | 0.0002 | (0.7742,0.8261) | 0.0368 | 0.0405 | (0.0937,0.8327) |
| 300 | -0.0003 | 0.0001 | (0.7782,0.8213) | -0.0072 | 0.0248 | (0.2022,0.7833) |
| 400 | -0.0004 | 0.0001 | (0.7809,0.8183) | -0.0033 | 0.0188 | (0.2403,0.7504) |
| 500 | -0.0007 | 0.0001 | (0.7825,0.8162) | -0.0094 | 0.0156 | (0.2618,0.7194) |
| $p=0.8, \beta=0.8$ | | | | | | |
| 50 | 0.0069 | 0.0005 | (0.7611,0.8526) | -0.1037 | 0.0773 | (0.2993,1.0932) |
| 100 | 0.0026 | 0.0002 | (0.7705,0.8348) | -0.0534 | 0.0359 | (0.4664,1.0269) |
| 200 | 0.0025 | 0.0001 | (0.7806,0.8244) | -0.0242 | 0.0148 | (0.5801,0.9715) |
| 300 | 0.0014 | 0.0001 | (0.7833,0.8195) | -0.0231 | 0.0077 | (0.6148,0.9390) |
| 400 | 0.0007 | 0.0001 | (0.7852,0.8162) | -0.0121 | 0.0054 | (0.6515,0.9243) |
| 500 | 0.0009 | 0.0000 | (0.7871,0.8147) | -0.0126 | 0.0043 | (0.6671,0.9077) |

From the findings presented in Table 1 it is observed that while there is no issue with the estimation of parameter $p$ with both bias and MSE decreasing with increase in sample size, and confidence interval shrinking towards true value of $p$, for $\beta$ the results are not good for small values of $p$. It is to be noted that for us $p$ is the main parameter where as $\beta$ is the additional nuisance parameter, so there is no reason to be disheartened. In case of confidence interval (CI) of $\beta$, the LCL<0 and UCL>1 should be treated as 0 and 1 for obvious reason. We report the original values for visualizing the improvement as sample size increases. Even the larger values of MSE for $\beta$ computed from the finite samples here indicates the stretching of the CIs beyond the boundaries $(0,1)$.

**Remark4:** This issue does not have impact on the data fitting because in both the examples the sample size was fairly large.

## 5. Testing of hypothesis

The geometric distribution is nested in the skewed geometric distribution for $\beta = 1$. Therefore for a given data set a natural question regarding the necessity of the additional parameter $\beta$ can be settled by testing the hypothesis $H_0: \beta = 1$ vs $H_1: \beta \neq 1$. That is, we test the hypothesis that the sampled observations are from the geometric distribution.

For this we adopt the Likelihood ratio (LR) test. The test statistic is

$$-2\log(LR) = -2\log L(\tilde{p},1|\mathrm{x})/L(\hat{p},\hat{\beta}|\mathrm{x}) \sim \chi^2 \text{ with 1 degree of freedom.}$$

where $\tilde{p} = \bar{x}/(1+\bar{x})$ is the mle of the parameter $p$ of geometric distribution in (1.2) and $\hat{p},\hat{\beta}$ are respectively the mles of $p,\beta$. If the calculated value of the test criterion exceeds 3.841 we reject the hypothesis at 5% level of significance.

The power curves of the LR test is given in Figure 3 for four different choices of the parameter $p$ with different sample sizes, considering $\beta \in$ {1.00, 0.85, 0.70, 0.55, 0.40, 0.25, 010}.

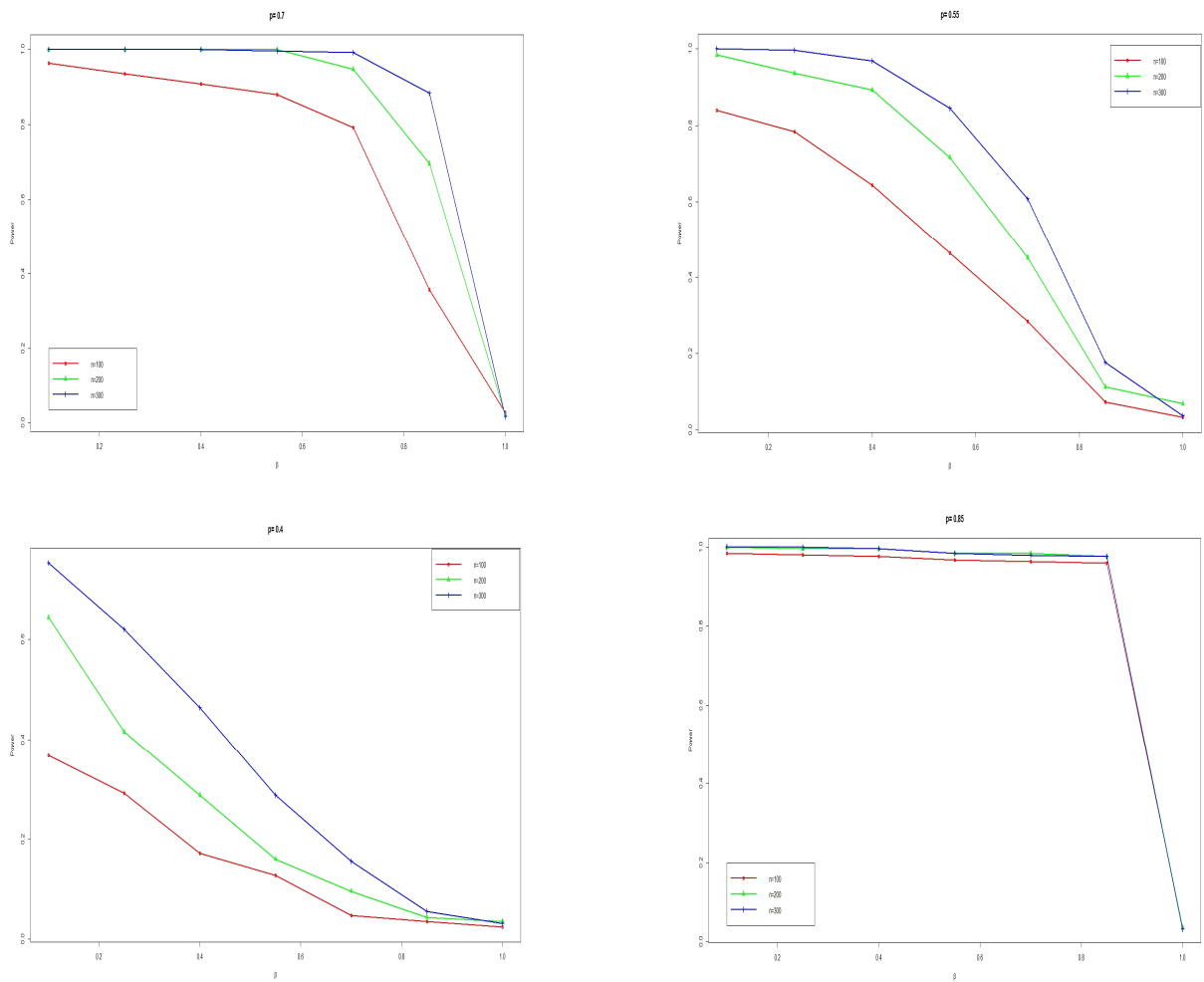

**Figure3**. Power curves for different sample sizes and choices of $p$.

**Observation from the power study:**

From the figures it is clear that LR test performs reasonably well for moderate to large values of $p$ for a fixed sample size as expected in line with our findings in performance of MLEs. Also as expected the power increases with sample size. We refrain for presenting the power curves for $p < 0.25$ and $p > 0.70$ as in former the performance with respect to power in very poor while in the later case the power for any $\beta < 1$ is almost equal to 1.

## 6. Data modelling applications

We consider the two datasets investigated by Bhati and Joshi (2018) for model fitting and compare $RSG(p,\beta)$ with the all the distributions considered by them including their weighted geometric (*WG*) distribution. The following list provides the competing distributions along with their pmf.

- $WG(p,\beta): P(X = x) = (1 - \beta)(1 - \beta^{p+1})\beta^x(1 - \beta^{p(x+1)})/(1 - \beta^p)$. (Bhati and Joshi, 2018)
- $NB(p,\beta): P(X = x) = \binom{p + x - 1}{x}\beta^x(1 - \beta^x)$ (Johnson, 2005)
- $ND(p,\beta): P(X = x) = (\log(1 - p\beta^x) - \log(1 - p\beta^{x+1}))/\log(1 - p)$ (Gomez, et al., 2011)
- $NGPL(p,\beta): P(X = x) = (\beta^2/(\beta + p)(1 + \beta)^{x+1})(1 + p(x + 1)/(1 + \beta))$ (Bhati et al., 2015)

The first data set originally from Klugman et al. (2012) concerns the number of claims made by automobile insurance policyholders, while the second dataset originally from Fisher (1941), is about frequency distribution of the number of ticks counted on 82 sheep. ML method of estimation is employed to estimate the model parameters. The findings are presented in Table 2 and Table 3.

**Table 2.** Number of claims made by automobile insurance policyholders

| # claims | Frequency | $RSG(p,\beta)$ | $WG(p,\beta)$ | $NB(p,\beta)$ | $ND(p,\beta)$ | $NGPL(p,\beta)$ |
|---|---|---|---|---|---|---|
| 0 | 1563 | 1564.687 | 1564.27 | 1564.54 | 1563.70 | 1564.57 |
| 1 | 271 | 265.32 | 265.12 | 264.58 | 266.15 | 264.28 |
| 2 | 32 | 38.48 | 39.05 | 39.44 | 38.75 | 39.69 |
| 3 | 7 | 5.57 | 5.62 | 5.66 | 5.50 | 5.59 |
| 4 | 2 | 0.96 | 0.94 | 0.78 | 0.90 | 0.87 |
|  | 1875 | 1875 | 1875 | 1875 | 1875 | 1875 |
| $(\chi_d^2, d)$ |  | (2.742, 2) | (2.938, 2) | (3.786, 2) | (3.002, 2) | (3.486, 2) |
| $p - val$ |  | 0.2538 | 0.230 | 0.151 | 0.223 | 0.176 |
| MLEs |  | (0.145, 0.001) | (0.873, 0.143) | (1.309, 0.871) | (-0.454, 0.141) | (8.835, 7.874) |
| $Var(\hat{p}), Var(\hat{\beta})$ |  | 0.0002, 1.618 | 0.683, 0.022 | 1.081, 0.214 | 0.495, 0.019 | 2.601, 0.439 |
| $Cov(p,\beta)$ |  | 0.0159 |  |  |  |  |
| AIC |  | 1990.58 | 1990.77 | 1991.00 | 1990.78 | 1991.18 |
| *Likelihood Ratio Test*: $H_0: \beta = 1$. That is sampled data comes from Geometric distribution with parameter $p$. The computed LR criterion is $-2\log likelihood\ ratio = \lambda = 1.2502$, is less than the cut-off value. Hence, $H_0$ Accepted. ||||||| 

**Table 3.** Number of ticks counted on 82 sheep

| # ticks | Frequency | $RSG(p,\beta)$ | $WG(p,\beta)$ | $NB(p,\beta)$ | $ND(p,\beta)$ | $NGPL(p,\beta)$ |
|---|---|---|---|---|---|---|
| 0 | 4 | 3.17 | 5.36 | 5.26 | 5.46 | 7.61 |
| 1 | 5 | 7.89 | 7.72 | 7.35 | 7.12 | 7.66 |
| 2 | 11 | 9.21 | 8.41 | 8.03 | 7.75 | 7.53 |
| 3 | 10 | 8.99 | 8.21 | 7.96 | 7.76 | 7.25 |
| 4 | 9 | 8.15 | 7.57 | 7.48 | 7.40 | 6.83 |
| 5 | 11 | 7.12 | 6.75 | 6.80 | 6.81 | 6.31 |
| 6 | 3 | 6.10 | 5.90 | 6.04 | 6.12 | 5.72 |
| 7 | 5 | 5.16 | 5.08 | 5.28 | 5.40 | 5.10 |
| 8,9,10 | 7 | 11.01 | 11.10 | 11.77 | 12.12 | 11.69 |
| 11,12,13,1,4 | 9 | 7.88 | 8.16 | 8.69 | 8.94 | 8.84 |
| ≥ 15 | 8 | 7.32 | 7.74 | 7.36 | 7.16 | 7.45 |
| | 82 | 82 | 82 | 82 | 82 | 82 |
| $(\chi_d^2, d)$ | | (7.2035, 8) | (8.476,8) | (9.124,8) | (9.844,8) | (12.666,8) |
| $p-val$ | | 0.5148 | 0.3884 | 0.3320 | 0.2761 | 0.1239 |
| MLEs | | (0.833, 0.601) | (0.834,1.759) | (1.777,0.271) | (1.276,0.311) | (2.312,0.808) |
| $Var(\hat{p}), Var(\hat{\beta})$ | | 0.0003, 0.0422 | 0.0008,2.5250 | 0.1211,0.0034 | 0.1600,0.0038 | 0.6839,0.0010 |
| $Cov(p,\beta)$ | | -0.0021 | | | | |
| AIC | | 477.92 | 478.98 | 479.92 | 480.88 | 483.44 |

*Likelihood Ratio Test*: $H_0: \beta = 1$. That is sampled data comes from Geometric distribution with parameter $p$. The computed LR criterion is $-2\log likelihood\ ratio = \lambda = 10.434$, is greater than the cut-off value. Hence $H_0$ is rejected and we conclude that the sampled observations are from $SG(p, \beta)$.

From the chi-square goodness-of-fit test reported in Table 2 and Table 3 it is obvious that the proposed $RSG(p, \beta)$. gives good fit to the data sets considered. Moreover in terms of the model selection criterion by AIC, SG is the best when compared with the other four extensions of geometric considered here.

## 7. Conclusion and discussion

Azzalani's skewing technique is used for the first time for generalizing a discrete distribution in this work, by considering the geometric distribution as baseline. The proposed model is very flexible and possesses compact forms for pmf, cdf, survival function, probability generating function and moments. It is an over dispersed distribution and can be long tailed, zero inflated, zero deflated for appropriate choices of the parameters. Problem of parameter estimation is discussed in detail and the performance is found to be satisfactory. Thus, it is envisaged that the proposed two-parameter model will find its application in count data modelling.

*The focus of this work was to introduce the skewing technique in the realm of discrete distributions and study the same with the geometric distribution.* Of course, the idea can be used for other count distributions as well. For the proposed model, other frequentist methods of parameter estimation may be attempted and Bayesian inference is also worth investigation when prior beliefs for the parameters are available. Modelling count data with covariates using some re-parametrized version of the proposed model will be taken up in follow-up work.